# Comment: Complex Causal Questions Require Careful Model Formulation: Discussion of Rubin on Experiments with "Censoring" Due to Death

**Stephen E. Fienberg**

## 1. INTRODUCTION

I am very pleased to be able to offer some reactions to yet another masterful paper by Rubin on the topic of causal inference. It is a special honor to do so because the paper was initially presented at Carnegie Mellon University as the 2005 Morris H. DeGroot Memorial Lecture and I was in the audience. "Morrie" DeGroot was my colleague, collaborator and close friend. He always raised questions about the naive use of randomization to answer causal questions from experiments. Thus, I believe that, had he been there to offer his own discussion at Rubin's oral presentation, he might have opined on the two issues that I address below, although perhaps with more wit.

The present paper fits quite naturally with Rubin [12], where he discusses the problematic nature of intermediate outcomes and R. A. Fisher's failure to recognize this problem. But this paper departs from that earlier one by avoiding the presentation of the key ideas using formal notation and equations. This makes for an interesting story but also for difficulties when one tries to follow the argument. The recipe for the resolution of most complex causal questions, we are told, is to frame them using potential outcomes and principal stratification. This is all well and good, but I still am not sure how to follow Rubin's recipe, either for the stylized example he uses or in other settings in the future.

*Stephen E. Fienberg is Maurice Falk University Professor of Statistics and Social Science, Department of Statistics, Carnegie Mellon University, Pittsburgh, Pennsylvania 15213, USA e-mail: fienberg@stat.cmu.edu.*

## 2. ALTERNATIVE REPRESENTATIONS FOR CAUSAL INFERENCE

I share Rubin's enthusiasm for representing causal questions using the formal framework of counterfactuals of which our philosophy colleagues are so fond. Rubin refers to these using the label "potential outcomes," harking back to Neyman's [7]. The reason I like this counterfactual representation is that it forces one to represent everything in terms of random variables, including randomization or any other allocation or missingness (censoring) mechanism; see [9, 10, 11]. Unfortunately, counterfactuals by their very nature lead us to condition on "unobservables" and thus they violate de Finetti's [3] dictum that conditional probabilities only make sense when we condition on actual observables, not simply potential ones. This is at least in part why Dawid [1, 2] has attempted to present a framework for causal inference similiar to Rubin's but which avoids the counterfactual representation. Lauritzen [6] has a related graphical model approach to this which he links to Pearl's [8] notion of "fixing" treatments or causes; see the similar ideas in [14].

My own preference is, as I suggested above, for representing every quantity under consideration using random variables, whether observed or unobserved, and then displaying these in graphical form using the standard methods for directed acyclic graphs. Thus, the act of randomization has a corresponding random variable and its introduction changes the graphical representation of the problem, often breaking the links between a treatment and an outcome variable; for example, see [4], as well as the more complete justification in [15]. This has the virtue of sidestepping Pearl's "unnatural" embellishments to the notation and representation of causal effects.

This is a very long preamble to a plea: If the arguments in the current paper are truly to hold sway, then: (1) They must have formal representation so we can see precisely where the assumptions fit in, and (2) We need to formulate them using the different causal representations, not simply the potential outcome framework.

## 3. ALTERNATIVE DEFINITIONS FOR CAUSAL EFFECT AND THEIR IMPLICATIONS

Rubin's original arguments for the role of randomization in experiments (e.g., see [10]) explicitly argued for a definition of *average causal effect* (ACE) based on a difference of expectations, and this suggests that the definition is "model-free" although the expectations are of course with respect to distributions that link to a model. I have always been troubled by the seeming arbitrariness of this representation. Why not the ratio or some other function of the expectations (cf. [1])?

In fact, it is relatively simple to see that the definition of ACE is intimately tied to linear models, and in recent work Sfer [13] and Fienberg and Sfer [5] have shown that tying the definition of causal effect to a formal parametric model resolves many of the seeming issues of bias associated with the effects of covariates in the nonlinear model setting. This is especially important for binary outcomes and for the modeling of some forms of survival.

I would therefore argue that we need to recast the principal stratification component of the present paper in a formal modeling context and then represent the censoring mechanisms in model-based terms as well. Then I believe we might really have a take-home lesson from the present paper on how to think about complex issues of causation with intermediate outcomes in the future.

## REFERENCES

[1] Dawid, A. P. (2000). Causal inference without counterfactuals (with discussion). *J. Amer. Statist. Assoc.* **95** 407–448. MR1803167

[2] Dawid, A. P. (2004). Probability, causality and the empirical world: A Bayes–de Finetti–Popper–Borel synthesis. *Statist. Sci.* **19** 44–57. MR2082146

[3] de Finetti, B. (1937). La prévision: ses lois logiques, ses sources subjectives. *Ann. Inst. H. Poincaré* **7** 1–68. MR1508036

[4] Fienberg, S. E., Glymour, C. and Spirtes, P. (1995). Discussion of "Causal diagrams for empirical research," by J. Pearl. *Biometrika* **82** 690–692.

[5] Fienberg, S. E. and Sfer, A. M. (2006). Randomization, models, and the estimation of causal effects. Submitted for publication.

[6] Lauritzen, S. (2004). Discussion of "Direct and indirect causal effects via potential outcomes," by D. B. Rubin. *Scand. J. Statist.* **31** 189–192. MR2066248

[7] Neyman, J. (1923). On the application of probability theory to agricultural experiments. Essay on principles. Section 9. *Ann. Agric. Sci.* **10** 1–51. (In Polish.) [Reprinted in English with discussion by T. Speed and D. B. Rubin in *Statist. Sci.* **5** (1990) 463–480.] MR1092986

[8] Pearl, J. (2000). *Causality: Models, Reasoning, and Inference.* Cambridge Univ. Press. MR1744773

[9] Rubin, D. B. (1976). Inference and missing data (with discussion). *Biometrika* **63** 581–592. MR0455196

[10] Rubin, D. B. (1978). Bayesian inference for causal effects: The role of randomization. *Ann. Statist.* **6** 34–58. MR0472152

[11] Rubin, D. B. (2004). Direct and indirect causal effects via potential outcomes (with discussion). *Scand. J. Statist.* **31** 161–170, 189–198. MR2066246

[12] Rubin, D. B. (2005). Causal inference using potential outcomes: Design, modeling, decisions. *J. Amer. Statist. Assoc.* **100** 322–331. MR2166071

[13] Sfer, A. M. (2005). Randomization and causality. Unpublished Ph.D. dissertation, Facultad de Ciencias Económicas, Univ. Nacional de Tucumán, San Miguel de Tucumán, Argentina. (In Spanish.)

[14] Spirtes, P., Glymour, C. and Scheines, R. (2000). *Causation, Prediction, and Search,* 2nd ed. MIT Press, Cambridge, MA. MR1815675

[15] Spirtes, P., Glymour, C., Scheines, R., Meek, C., Fienberg, S. E. and Slate, E. (1999). Prediction and experimental design with graphical causal models. In *Computation, Causation, and Discovery* (C. Glymour and G. F. Cooper, eds.) 65–93. AAAI Press/MIT Press, Cambridge, MA. MR1689956